\documentclass[11pt]{amsart}

\usepackage[usenames]{color}

\usepackage[colorlinks, linkcolor=reference,
            citecolor=citation,
          urlcolor=e-mail]{hyperref}

\definecolor{e-mail}{rgb}{0,.40,.80}
\definecolor{reference}{rgb}{.20,.60,.22}
\definecolor{citation}{rgb}{0,.40,.80}
\usepackage{amsmath,times,pslatex,enumitem}
\usepackage{amsthm}
\usepackage{amsfonts}
\usepackage{amssymb}
\usepackage{latexsym}

\usepackage{hyperref}

\newtheorem{theorem}{Theorem}

\newtheorem{remark}{Remark}
\newtheorem*{remark*}{Remark}

\author[]{Mariya Bessonov}
\address{Department of Mathematics, New York City College of Technology, CUNY} \email{mbessonov@citytech.cuny.edu}\thanks{Research of the
first author was partially supported by the NSF grant DMS-1515800}

\author[]{Dima Grigoriev}
\address{CNRS, Math\'ematiques, Universit\'e de Lille, 59655, Villeneuve d'Ascq, France}
\email{dmitry.grigoryev@math.univ-lille1.fr}

\author[]{Vladimir Shpilrain}
\address{Department of Mathematics, The City College of New York, New York,
NY 10031} \email{shpil@groups.sci.ccny.cuny.edu}
\thanks{Research of the third author was partially supported by
the ONR (Office of Naval Research) grant N000141512164}

\begin{document}

\title{Probabilistic solution of Yao's millionaires' problem}

\begin{abstract}
We offer a probabilistic solution of Yao's millionaires' problem
that gives correct answer with probability (slightly) less than 1
but on the positive side, this solution does not use any one-way
functions.

\end{abstract}

\maketitle

\section{Introduction}

The ``two millionaires problem" introduced by Yao in \cite{Yao} is:
\begin{quote}
Alice has a private number  $a$ and Bob has a private number $b$,
and the goal of the two parties is to solve the inequality $a  \le b
?$ without revealing the actual values of $a$ or $b$, or more
stringently, without revealing any information about $a$ or $b$
other than $a  \le b$ or $a > b$.\end{quote} We note that all known
solutions of this problem (including Yao's original solution) use
one-way functions one way or another. Informally, a function is {\it
one-way} if it is efficient to compute but computationally
infeasible to invert on ``most" inputs. One problem with those
solutions is that it is  not known whether one-way functions
actually exist, i.e., the functions used in the aforementioned
solutions are just {\it assumed} to be one-way. Also, solutions that
use one-way functions inevitably use assumptions of limited
computational power of the parties, and this assumption is arguably
more ``physical" than ``mathematical" in nature, although there is a
mathematical theory of computational complexity with a (somewhat
arbitrary) focus on distinction between polynomial-time and
superpolynomial-time complexity of algorithms.

Speaking of physics, in our earlier papers \cite{Yao_physics} and
\cite{GKS} we offered several solutions of Yao's millionaires'
problem without using one-way functions, but 
using real-life procedures (not implementable on a Turing machine).
What is important is that some of these solutions can be used to
build a public-key encryption protocol secure against a
computationally unbounded (passive) adversary, see \cite{GKS}.

Here we make an assumption that both private numbers $a$ and $b$ are
uniformly distributed on integers in a public interval $[1, n]$.
This assumption may be questionable as far as the millionaires'
problem itself is concerned, but we keep in mind potential
applications to cryptographic primitives, in which case the above
assumption could be just fine. We also note that this assumption can
be relaxed to $a$ and $b$ being identically (not necessarily
uniformly) distributed on integers in $[1, n]$ because in that case,
a monotone function ${\mathcal F}$ can be applied to both $a$ and $b$ so that
${\mathcal F}(a)$ and ${\mathcal F}(b)$ become uniformly distributed (on a different
interval though), see \cite[Section 2.2.1]{Glass}. This is called
the {\it inverse transform method.}

We also note that our solution of the millionaires' problem has a
``symmetric" as well as ``asymmetric" version. In the ``asymmetric"
version, only one of the two parties ends up knowing whether $a < b$
or not, even if the other party is computationally unbounded. (Of
course, she can then share this information with the other party if
she chooses to.) This implies that a third party (a passive
observer) will not know whether $a < b$ either, and this is the key
property for building a public-key encryption protocol secure
against  a computationally unbounded (passive) adversary, see
\cite{GKS} for details.

The way our solution in this paper works is roughly as follows.
Alice applies a randomized function $F$ to her private number $a$
and obtains the result $A=F(a)$ that she either keeps private
(``asymmetric" version) or makes public (``symmetric" version). Bob
applies a randomized function $G$ to his private number $b$ and
obtains the result $B=G(b)$ that he makes public. Then Alice, based
on $a, A$ and $B$, makes a judgement whether $a < b$ or not.
Specifically, in our protocol in Section \ref{Protocol}, she
concludes that $a < b$ if and only if $A < B$. 
We show in Section \ref{Probabilities} that, with appropriate choice
of parameters, this judgement is correct with probability converging
to 1 as $n$ (the interval length) goes to infinity. Computer
simulations suggest (see Section \ref{parameters}) that this
convergence is actually rather fast.

\section{Protocol}
\label{Protocol}

Recall that in a {\it simple symmetric random walk}, a point on a
horizontal line moves one unit left with probability $\frac{1}{2}$
or one unit right with probability $\frac{1}{2}$. Below is our
protocol for a probabilistic solution of Yao's millionaires'
problem, under the assumption that both private numbers $a$ and $b$
are uniformly distributed on integers in a public interval $[1, n]$.

\begin{enumerate}[leftmargin=1.2cm,labelsep=-0.5cm,itemsep=0.2cm,align=left]
\item[1.] Alice's private number $a$ is the  starting
point of her random walk. Alice does a simple symmetric random walk
with $f(n)$ steps, starting at $a$. Let $A$ be the end point of
Alice's random walk. Alice either keeps $A$ private (``asymmetric"
version) or makes it public (``symmetric" version).

\item[2.] Bob's private number $b$ is the  starting
point of his random walk.  Bob  does a simple symmetric random walk
with $g(n)$ steps, starting at $b$. Let $B$ be the end point of
Bob's random walk. Bob makes $B$ public.

\item[3.] Alice concludes that $a < b$ if
and only if $A < B$.

\item[4.] In case Alice has published her $A$ (``symmetric" version), Bob, too, concludes that $a < b$ if
and only if $A < B$.

\end{enumerate}

We emphasize that in the ``asymmetric" version, only Alice ends up
knowing (with significant probability) whether or not $a < b$.
Neither Bob nor a third party observer end up knowing this
information unless Alice chooses to share it.

In this paper, we focus on the arrangement where $f(n)=g(n)$, i.e.,
the parties do the same number of steps. Other arrangements are
possible, too; in particular, as we remark in Section \ref{nowalk},
Alice can (slightly) improve the probability of coming to the
correct conclusion on $a  < b ?$ if she does not walk at all, i.e.,
if $f(n)=0$. This arrangement is ``highly asymmetric" but it is
useful to keep in mind for future work.

\section{Probabilities}
\label{Probabilities}

We start with the following

\begin{remark}If $a$ and $b$ are independent random
variables and each is uniformly distributed on $\{1,2, \ldots ,
n\}$, then the expected  value of $|a-b|$ is $E(|a-b|) = \frac{(n^2
- 1)}{3n}$, which is asymptotically equal to $\frac{n}{3}$.
\end{remark}

Indeed, note that $|a-b| = \max (a, b) - \min (a, b)$. By symmetry,
$$E(\max (a, b)) = n+1 - E(\min (a, b)),$$ hence $$E(|a-b|) = n+1 -
2E(\min (a, b)).$$ Now direct computation gives $$ E(\min (a, b)) =
\sum_{k=1}^n k \cdot (\frac{2}{n}\cdot \frac{n-k}{n} +
\frac{1}{n^2}) = n+1 + \frac{(n+1)(1-4n)}{6n}.$$ Then $$E(|a-b|) = n+1
- 2E(\min (a, b)) = \frac{(n^2 - 1)}{3n}.$$

\begin{remark}
If $t$ is a positive integer and $a$ and $b$ are independent and
uniformly distributed on $\{1,2, \ldots, n\}$, then $P(|a-b| < t) <
2t/n$.
\end{remark}
To see this,
\begin{align*}
P(|a-b| < t) & = \sum_{j=1}^{n}P\left(\{|a-b| < t\} \cap \{a = j\}\right) \\
& < \sum_{j=1}^{n} P\left(b\in\{j-(t-1),j-t+2,\ldots, j+t-1\} \cap a = j\right)\\
& = \sum_{j=1}^{n} P\left(b\in\{j-t+1,j-t+2,\ldots, j+t-1\}\right)P(a = j)\\
&= \frac{2t}{n}\cdot\frac{1}{n}\cdot n = \frac{2t}{n}
\end{align*}
Thus, $a$ and $b$ are likely to be sufficiently far apart, which
explains why the probability that our solution is correct is
sufficiently high.

Recall that the probability of our solution being correct is the
conditional probability $P(a < b ~| ~A < B)$. It depends on the
functions $f(n)$ and  $g(n)$, and we consider a couple of cases
here, focusing on the case where $f(n)=g(n)$. First we recall

\begin{theorem}[{see \cite[Theorem~2.2 and Remark~2.6]{EL03}}]
Let $S_m$ be the location of the simple symmetric random walk on
$\mathbb{Z}$ after $m$ steps with $S_0 = 0$. Let $1/2 < \alpha < 1$.
Then
\begin{align}
\lim_{m\to\infty} m^{1-2\alpha}\log P\left(S_m > xm^{\alpha}\right)
= -\frac{x^2}{2}
\end{align}
\end{theorem}

That is, for large enough $m$,
$$P\left(S_m > xm^{\alpha}\right) \approx e^{-\frac{x^2m^{2\alpha - 1}}{2}} \to 0 \text{ as } m\to\infty,$$
for $x=1$,
$$P\left(S_m > m^{\alpha}\right) \approx e^{-\frac{m^{2\alpha - 1}}{2}} \to 0 \text{ as } m\to\infty,$$
so the probability that the displacement from the starting point is
greater than $O(\sqrt{m})$ tends to $0$.

If $g(n) = O(n^{2-2\epsilon})$, $n$ is fixed, and $\epsilon\to 0$,
then $P(a<b | A<B)$ will not be close to $1$ since the typical
displacement is $O(n^{1-\epsilon})$. This probability $P(a<b | A<B)$
will tend to $1$ for any fixed $\epsilon > 0$ and $n\to\infty$.

For $m$ fixed and $\alpha\in[1/2,1]$, one has
\begin{align*}
P(S_m \geq m^{\alpha}) & = P\left(e^{S_mm^{\alpha-1}} \geq e^{m^{2\alpha-1}}\right)\\
&\leq  \exp{\left(-m^{2\alpha -1}\right)}E\left(\exp{\left(S_mm^{\alpha -1}\right)}\right)\\
&=  \exp{\left(-m^{2\alpha -1}\right)}\left(E\left(e^{\left(Xm^{\alpha -1}\right)}\right)\right)^m\\
&= \exp{\left(-m^{2\alpha -1} + m\ln\left(\cosh(m^{\alpha -
1})\right)\right)} \leq \exp{\left(-\frac{m^{2\alpha -
1}}{2}\right)},
\end{align*}
where $X$ is a random variable taking on $1$ and $-1$ with equal
probability (the step distribution of the simple symmetric random
walk). The first inequality is an application of Markov's inequality
and the last inequality holds because $\ln\cosh(x) \leq
\frac{x^2}{2}$ for all $x\in\mathbb{R}$.

It will follow that if there are $m=n^{\lambda}$ steps in
Alice's and in Bob's random walks, then for any
$$\alpha\in(1/2,\min\{1,\ln\left(n/2\right)/\lambda\ln(n)\}),$$ one
has
\begin{align}\label{eq:main}P&(a<b | A<B)\notag\\ &\geq \left(1 - \exp{\left(-\frac{m^{2\alpha - 1}}{2}\right)}\right)^2\left(\frac{(n-2m^{\alpha}+1)(n-2m^{\alpha})}{n^2-n}\right)\left(1 - \frac{1}{n}\right),
\end{align}
which approaches $1$ in the limit as $n\to\infty$.

To see why \eqref{eq:main} is true, consider
\begin{align}\label{eq:condProb}
P(a<b | A < B) = \frac{P\left(\{a < b\}\cap\{A < B\}\right)}{P(A <
B)} = \frac{P(A<B | a<b)P(a<b)}{P(A < B)}.
\end{align}
The denominator of \eqref{eq:condProb}, $P(A < B) < 1/2$. Indeed,
given that Alice's and Bob's random walks have the same number of
steps and are denoted $A_k$ and $B_k$, the difference between their
random walks $Y_k = B_k - A_k$ for $k=0,1,\ldots,m$ is a lazy
symmetric random walk with probability $1/2$ of staying in place and
probabilities $1/4$ each of moving two steps to the left or to the
right. Since the random walks are symmetric and the starting points
$a$ and $b$ are selected uniformly and independently of each other
on the same interval, by symmetry, $$P(A < B) = P(Y_m > 0) = P(Y_m <
0) = P(A > B),$$ and $P(A=B) = P(Y_m = 0) = O(1/\sqrt{m})$.

Then, $P(a<b) = \frac{n^2 - n}{2n^2}$ since $a$ and $b$ are chosen
independently and uniformly at random on $\{1,2,\ldots,n\}$. There
are $n^2$ different ordered pairs $(a,b)$, n of which have $a=b$,
and half of the remaining $n^2-n$ ordered pairs have $a<b$.

For the other term in the numerator of \eqref{eq:condProb}, let $E$
be the event that $b - a \geq 2m^{\alpha}$ and assume $\alpha$ is
such that $2m^{\alpha}$ is an integer. Then
\begin{align}\label{eq:condSwitch}
P(A < B | a < b) &\geq P(\{A < B\} \cap E | a < b)\\ \notag & =
P\left(A < B | \{b-a \geq 2m^{\alpha}\}\cap\{a<b\}\right)\cdot
P\left(b - a \geq 2m^{\alpha} | a < b\right)\\ \notag & = P\left(A <
B | b-a \geq 2m^{\alpha}\right)\cdot P\left(b - a \geq 2m^{\alpha} |
a < b\right).
\end{align}

For the second term in \eqref{eq:condSwitch}, recall that there are
$\frac{n^2-n}{2}$ ordered pairs $(a,b)$ with $a<b$, each ordered
pair equally likely. Thus, we have
\begin{align*}
P\left(b-a \geq 2m^{\alpha} | a < b\right) &= \sum_{j=2m^{\alpha}}^{n-1}P\left(b-a = j | a < b\right)\\
&=\frac{(n-2m^{\alpha}) + (n-2m^{\alpha}-1)+ \ldots + 1}{\frac{n^2-n}{2}}\\
& = \frac{2\cdot\frac{(n-2m^{\alpha}+1)(n-2m^{\alpha})}{2}}{n^2-n}\\
& = \frac{(n-2m^{\alpha}+1)(n-2m^{\alpha})}{n^2-n}.
\end{align*}
If $m=n^{\beta}$, this expression is greater than $0$ when $\alpha <
\frac{\ln \left(n/2\right)}{\lambda\ln n}$.

Let $F_{\alpha}$ be the event that each of Alice's and Bob's random
walks traveled distance no more than $m^{\alpha}$. Each random walk
has probability less than $e^{-m^{2\alpha - 1}/2}$ of traveling more
than $m^{\alpha}$ from its starting point. Thus, since the distance
traveled of each walk is independent of the starting points and
since the random walks are independent of each other,
\begin{align*}
P\left(A < B | b-a \geq 2m^{\alpha}\right) &\geq P(F_{\alpha} | b-a \geq 2m^{\alpha}) = P(F_{\alpha}) \\
&\geq \left(1 -\exp{\left(-\frac{m^{2\alpha-1}}{2}\right)}\right)^2.
\end{align*}
Then \eqref{eq:main} follows.

An improvement to the lower bound on $P(a<b | A <B)$ for smaller
values of $n$ can be obtained by improving the bound in
\eqref{eq:condSwitch}
\begin{align}\label{eq:condSwitch2}
P(A < B | a < b) &= P(\{A < B\} \cap E | a < b) + P(\{A<B\}\cap E^c
| a < b)
\end{align}
For the second term on the right of \eqref{eq:condSwitch2},
\begin{align}\label{eq:complementBound}
P(\{A<B\}\cap E^c | a < b)& = P\left(A < B | 0 < b-a <
2m^{\alpha}\right)\cdot P\left(b - a < 2m^{\alpha} | a < b\right)\\
\notag &
> (1/2)\cdot\left(1-\frac{(n-2m^{\alpha}+1)(n-2m^{\alpha})}{n^2-n}\right),
\end{align}
where we use $$P\left(A < B | 0 < b-a < 2m^{\alpha}\right) = P(Y_m >
0 | 0 < Y_0 < 2m^{\alpha}) > 1/2,$$  since $Y_k$ is a symmetric
random walk. Then, combining \eqref{eq:condProb},
\eqref{eq:condSwitch}, \eqref{eq:condSwitch2}, and
\eqref{eq:complementBound},
\begin{align}\label{eq:main2}
P&(a<b | A<B)\notag\\ &\geq \left(1 - \exp{\left(-\frac{m^{2\alpha -
1}}{2}\right)}\right)^2\left(\frac{(n-2m^{\alpha}+1)(n-2m^{\alpha})}{n^2-n}\right)\left(1
- \frac{1}{n}\right)\\ \notag & + (1/2)\cdot\left(1 -
\frac{(n-2m^{\alpha}+1)(n-2m^{\alpha})}{n^2-n}\right)\left(1 -
\frac{1}{n}\right).
\end{align}

\subsection{What if Alice does not walk?} \label{nowalk}
If Alice  does not walk and makes a judgement based on her point $a$ and the terminal point $B$ of Bob's walk, then the probability in question is $P(a<b | a < B)$. The (somewhat informal) argument below shows that this probability is, in fact, greater than $P(a<b | A < B)$, although computer simulation shows that the difference is rather small. Thus, for the purpose of solving the millionaires problem itself, it is preferable to use $f(n)=0$  and $g(n)=n^{\frac{4}{3}}$ for the number of steps in Alice's and Bob's random walk, respectively. However, if one has in mind a possible conversion of such a solution to an encryption scheme, then having $f(n)=0$ is not optimal from the security point of view. We leave this discussion to another paper though, while here we explain why $P(a<b | a < B) > P(a<b | A < B)$.

Note that if Alice does not walk, the difference between Bob's and Alice's position is a simple symmetric random walk, $X_k$, with probability 1/2 each of moving to the right or left 1 step.

On the other hand,  if both walk, then the difference between Bob's and Alice's position is a lazy random walk, $Y_k$, with probability 1/2 of staying in place and 1/4 each of moving to the right or left 2 steps.
Then $$P(a<b | a < B) > P(a<b | A < B)\iff P(X_0 > 0 | X_m > 0) > P(Y_0 > 0 | Y_m > 0).$$

It is well known that the mean squared displacement is greater for the lazy walk. Specifically, in our situation,
$$E(X_m^2) = m,\ \ E(Y_m^2) = 2m.$$
Based on this, we indeed have $P(X_0 > 0 | X_m > 0) > P(Y_0 > 0 | Y_m > 0)$, so
$$P(a<b | a < B) > P(a<b | A < B).$$

\subsection{The case of $n^{\frac{4}{3}}$ steps in random walks} \label{43}
Using {\sc Maple}, we found the maximum over $\alpha\in(1/2,\frac{\ln
\left(n/2\right)}{\lambda\ln n})$ of \eqref{eq:main2}  for several values of $n$, with
$n^{4/3}$ steps in both Alice's and Bob's random walks (see Table~\ref{tab:1}).
\begin{table}[h!]
\centering
\begin{tabular}{|c|c|c|}
\hline
$n$ &$P(a<b | A<B) $ & $\alpha $ \\
\hline
$10^3$   & $\geq 0.586$ & $\approx 0.574$ \\
\hline
$10^4$&  $ \geq 0.743$& $ \approx 0.574$\\
\hline
$10^5$&  $\geq 0.859$& $ \approx 0.568$\\
\hline
$10^6$& $\geq 0.927$& $ \approx 0.563$\\
\hline
$10^7$&  $\geq 0.963$& $\approx 0.557$\\
\hline
$10^8$&  $\geq 0.982$& $\approx 0.553$\\
\hline
 $10^9$&  $ \geq 0.991$& $\approx 0.549$\\
 \hline
\end{tabular}
\caption{Case of $n^{4/3}$ steps}\label{tab:1}
\end{table}
We emphasize that these are lower bounds; the actual speed
of convergence to 1 appears to be faster. For example, computer
simulations suggest that already for $n=1000$,  ~$P(a<b | A<B)$ is
about $0.9$. If $n=2000$, then $P(a<b | A<B)$ is about $0.99$.

\subsection{The case of $n^{\frac{5}{3}}$  steps in random walks}
The maximum over $\alpha\in(1/2,\frac{\ln \left(n/2\right)}{\lambda\ln
n})$ of \eqref{eq:main2} for several values of $n$ in this case is given in Table~\ref{tab:2}.
\begin{table}[h!]
\centering
\begin{tabular}{|c|c|c|}
\hline
$n$ &$P(a<b | A<B) $ & $\alpha $ \\
\hline
$10^3$   & $ \geq 0.453$& $ \approx 0.500$ \\
\hline
$10^4$&  $ \geq 0.466$& $\approx 0.517$\\
\hline
$10^5$&  $\geq 0.514$& $\approx 0.526$\\
\hline
$10^6$& $\geq 0.586$& $\approx 0.529$\\
\hline
$10^7$&  $ \geq 0.667$& $\approx 0.530$\\
\hline
$10^8$&  $\geq 0.743$& $\approx 0.530$\\
\hline
 $10^9$&  $\geq 0.807$& $\approx 0.529$\\
 \hline
\end{tabular}
\caption{Case of $n^{5/3}$ steps}\label{tab:2}
\end{table}
Again, the actual speed of convergence to 1 appears to be
faster. Computer simulations suggest that for $n=1000$,  ~$P(a<b |
A<B)$ is about $0.75$ in this case.


\subsection{The probability to guess the other party's
number}\label{guess}
Another probability that we are interested in
is the probability for Alice to correctly guess Bob's private number
$b$ based on the public $B$. The most likely position of the point
$b$ is $b=B$ (assuming that $g(n)$ is even), and the probability for
that to actually happen is (using Stirling's formula) approximately
$\sqrt{\frac{2}{\pi g(n)}}$. Thus, we have:
\begin{enumerate}[leftmargin=1.2cm,labelsep=-0.5cm,itemsep=0.2cm,align=left]
\item[1.] For $f(n) = g(n) = n$, Alice's best guess for $b$ has
probability about $\sqrt{\frac{2}{\pi n}}$ to be correct.

\item[2.] For $f(n) = g(n) = n^{\frac{4}{3}}$, Alice's best guess for $b$ has
probability about $\sqrt{\frac{2}{\pi n^{\frac{4}{3}}}} =
\frac{\sqrt{2}}{n^{\frac{2}{3}} \sqrt{\pi} }$ to be correct.

\item[3.] For $f(n) = g(n) = n^{\frac{5}{3}}$, Alice's best guess for $b$ has
probability about $\sqrt{\frac{2}{\pi n^{\frac{5}{3}}}}$ to be correct.

\end{enumerate}

These probabilities can be compared to the {\it a priori}
probability for either party to guess the other party's number
correctly (with or without knowing the probability distribution),
which is
\begin{equation}\label{eq:probapprox}
\frac{1}{n}\sum_{k=1}^n \frac{1}{k} \approx \frac{\ln
n}{n}.
\end{equation} Indeed, if the range for $a$ and $b$ is $[N_1, N_2]$, and
Bob's integer $b$ happens to be equal to $N_1$, then, after having
found out that $a \le  b$, Bob knows that Alice's integer is $a =
N_1$. Then, if $b=N_1+1$, the information $a \le  b$ tells Bob that
either $a = N_1$  or  $a = N_1+1$, so he can guess $a$ correctly
with probability 1/2. Thus, 
 in the ``ideal" situation where an oracle just tells
Bob that, say, $a \le b$, the total probability for Bob to guess $a$
correctly is~\eqref{eq:probapprox}.

As another point of comparison, we mention a very simple solution of
the millionaires' problem from \cite{GKS}:

\begin{enumerate}[leftmargin=1.2cm,labelsep=-0.5cm,itemsep=0.2cm,align=left]

\item Alice begins by breaking the set of $n$ integers from the
interval $[1, n]$ into approximately $\sqrt{n}$ subintervals with
approximately $\sqrt{n}$ integers in each, in such a way that her
integer $a$ is an endpoint of one of the subintervals.

\item Alice then sends the  endpoints of all the subintervals  to
Bob. (Alternatively, she can send just a compact description of the
endpoints.)

\item Bob tells  Alice in which subinterval his integer $b$ is. By the above
property of Alice's subintervals, all elements of the subinterval
pointed at by Bob are either less than (or equal to) $a$ or greater
than $a$, so Alice now has a solution of the  inequality  $a \le b
?$.

\end{enumerate}

It is obvious that the probability for Bob to guess Alice's integer
$a$ correctly, as well as the probability for Alice to guess Bob's
integer $b$ correctly, is approximately $\frac{1}{\sqrt{n}}$.

As a side remark, we note that in this solution Alice ends up with
exactly the same information about Bob's number $b$ as a third party
observer does, and this information is deterministic, so Alice does
not get any advantage over a third party in case she is thinking of
using this solution to send encrypted information to Bob. See
\cite{GKS} for details on situations where a solution of the
millionaires' problem can be used to build a public-key encryption
scheme.


\section{Suggested parameters for practical use and experimental results}
\label{parameters}

We recommend selecting an interval of length $n = 8000$ and
selecting $n^{\frac{4}{3}} = 160,000$  steps in the parties' random
walks. If  $a$ and $b$ are uniformly distributed on integers in an
interval $[1, N]$ with $N<n$, then they are identically (although
not uniformly) distributed on integers in $[1, n]$, in which case
one can use the {\it inverse transform method} mentioned in our
Introduction to reduce to the case of the uniform distribution on
$[1, n]$. If $N>n$, then the parties can represent their private
numbers in the form $\sum c_k n^k$ with $c_k < n$ and compare the
coefficients $c_k$, starting with the largest $k$.

With $n = 8000$ and $m = 160,000$ steps, the probability for Alice to
guess Bob's private number $b$ is $\frac{\sqrt{2}}{n^{\frac{2}{3}}
\sqrt{\pi} } \approx 0.002$ (see our Section \ref{guess}) and,
according to computer simulations, $P(a<b | A<B) \approx 0.99$.

With these parameters, simulation of a random walk takes 0.05 sec on
a regular desktop computer.

\section{Conclusions}

Recall that $n$ is the length of an interval from which the two
parties' private integers are selected.

\begin{itemize}[leftmargin=1.2cm,labelsep=-0.5cm,itemsep=0.2cm,align=left]
\item We see that, when choosing $n^{\lambda}$ steps of the
parties' random walks, $\lambda$ should be less than 2 for $P(a<b |
A<B)$ to converge to 1 as $n \to \infty$. If $\lambda \ge 2$, then
$P(a<b | A<B)$ does not converge to 1 as $n\to\infty$.

\item In choosing a particular $\lambda < 2$, there is a
trade-off between the probability for Alice to correctly solve $a
\le b ?$ and the the probability for Alice to guess Bob's private
number $b$. More specifically, the closer the number of steps of the
parties' random walks is to $n^2$, the slower is the convergence of
$P(a<b | A<B)$ to 1, but at the same time, the bigger the spread of
the public point $B$ around the private point $b$ is, thus reducing
the probability for Alice to guess Bob's number $b$.

\item We choose $n^{\frac{4}{3}}$ steps as the ``equilibrium" in
this trade-off. Lower bounds for $P(a<b | A<B)$ in this case, as
computed in our Section \ref{43}, can make an impression that our
method is very inefficient since $n$ has to be very large for $P(a<b
| A<B)$ to become close to 1. However, the actual speed of
convergence to 1 appears to be faster. For example, computer
simulations suggest that already for $n=1000$,  ~$P(a<b | A<B)$ is
about $0.9$ in that case.

\item Our recommendation for the choice of parameters is: $n =
8000$, and the number of steps in the parties' random walks is
$n^{\frac{4}{3}} = 160,000$. With these parameters, $P(a<b | A<B)
\approx 0.99$, and the probability for Alice to guess Bob's private
number $b$ is about $\frac{\sqrt{2}}{n^{\frac{2}{3}} \sqrt{\pi} } \approx
0.002$.

\item In this paper, the focus is on the arrangement where Alice and
Bob do the same number of steps in their random walks. Other
arrangements are possible, too; in particular, as we remark in
Section \ref{nowalk}, Alice can (slightly) improve the probability
of coming to the correct conclusion on $a  \le b ?$ if she does not
walk at all. This arrangement is ``highly asymmetric" but it is
nevertheless useful to keep in mind.

\end{itemize}

\end{document}